\documentclass[12pt]{amsart}

\usepackage{amsmath,amssymb,amsthm,fullpage,verbatim}
\usepackage{latexsym}
\usepackage{pstricks}

\newcommand{\ES}{\text{ES}}

\newcommand{\quat}[2]{\displaystyle{\biggl(\frac{#1}{#2}\biggr)}}
\newcommand{\legen}[2]{\left(\frac{#1}{#2}\right)}

\DeclareMathOperator{\MM}{M}
\DeclareMathOperator{\ord}{ord}

\DeclareMathOperator{\Hom}{Hom}

\DeclareMathOperator{\Gal}{Gal}

\DeclareMathOperator{\Div}{Div}

\DeclareMathOperator{\GL}{GL}

\DeclareMathOperator{\SL}{SL}

\DeclareMathOperator{\PGL}{PGL}

\DeclareMathOperator{\tr}{tr}

\DeclareMathOperator{\Coind}{Coind}

\DeclareMathOperator{\M}{M}
\DeclareMathOperator{\nrd}{nrd}

\DeclareMathOperator{\Frob}{Frob}

\newcommand{\cp}{\times}
\newcommand{\tp}{\otimes}
\newcommand{\TP}{\bigotimes}
\newcommand{\ds}{\oplus}

\newcommand{\bs}{\backslash}

\newcommand{\subs}{\subset}

\newcommand{\ZZ}{\mathbb{Z}}
\newcommand{\QQ}{\mathbb{Q}}
\newcommand{\CC}{\mathbb{C}}
\newcommand{\RR}{\mathbb{R}}
\newcommand{\PP}{\mathbb{P}}
\newcommand{\FF}{\mathbb{F}}

\newcommand{\Z}{\mathbb{Z}}
\newcommand{\Q}{\mathbb{Q}}
\newcommand{\C}{\mathbb{C}}
\newcommand{\R}{\mathbb{R}}
\newcommand{\F}{\mathbb{F}}

\newcommand{\cF}{\mathcal{F}}
\newcommand{\calF}{\mathcal{F}}

\newcommand{\cH}{\mathcal{H}}

\newcommand{\cO}{\mathcal{O}}
\newcommand{\calO}{\mathcal{O}}

\newcommand{\cS}{\mathcal{S}}

\newcommand{\fa}{\mathfrak{a}}
\newcommand{\fA}{\mathfrak{A}}
\newcommand{\fb}{\mathfrak{b}}

\newcommand{\fD}{\mathfrak{D}}

\newcommand{\fM}{\mathfrak{M}}
\newcommand{\fn}{\mathfrak{n}}
\newcommand{\fN}{\mathfrak{N}}

\newcommand{\fp}{\mathfrak{p}}
\newcommand{\fP}{\mathfrak{P}}
\newcommand{\fq}{\mathfrak{q}}

\newcommand{\fr}{\mathfrak{r}}

\newcommand{\ga}{\alpha}

\renewcommand{\gg}{\gamma}
\newcommand{\gG}{\Gamma}

\newcommand{\gi}{\iota}

\newcommand{\gm}{\mu}

\newcommand{\go}{\omega}

\newcommand{\gp}{\pi}

\newcommand{\gt}{\tau}

\numberwithin{equation}{section}

\newtheorem{prop}[equation]{Proposition}

\newtheorem{theorem}[equation]{Theorem}

\newtheorem*{theorema*}{Theorem}
\newtheorem*{theoremb*}{Theorem B}

\theoremstyle{remark}

\newtheorem{remark}[equation]{Remark}

\theoremstyle{definition}

\newtheorem{alg}[equation]{Algorithm}

\newcommand{\slsh}[1]{\,|_{#1}\,}

\newenvironment{enumalg}
{\begin{enumerate}}
{\end{enumerate}}

\begin{document}
\title[Computing Hilbert modular forms]{Computing systems of Hecke eigenvalues associated to Hilbert modular forms}

\author{Matthew Greenberg}
\address{University of Calgary, 2500 University Drive NW, Calgary, AB, T2N 1N4, Canada}
\email{mgreenbe@math.ucalgary.ca}
\author{John Voight}
\address{Department of Mathematics and Statistics, University of Vermont, 16 Colchester Ave, Burlington, VT 05401, USA}
\email{jvoight@gmail.com}
\date{\today}

\begin{abstract}
We utilize effective algorithms for computing in the cohomology of a Shimura curve together with the Jacquet-Langlands correspondence to compute systems of Hecke eigenvalues associated to Hilbert modular forms over a totally real field $F$.
\end{abstract}

\maketitle

The design of algorithms for the enumeration of automorphic forms has emerged as a major theme in computational arithmetic geometry.  Extensive computations have been carried out for elliptic modular forms and now large databases exist of such forms \cite{CremonaLatest,SteinWatkins}.  As a consequence of the modularity theorem of Wiles and others, these tables enumerate all isogeny classes of elliptic curves over $\QQ$ up to a very large conductor.  The algorithms employed to list such forms rely heavily on the formalism of modular symbols, introduced by Manin \cite{Manin} and extensively developed by Cremona \cite{Cremona}, Stein \cite{Stein}, and others.  For a positive integer $N$, the space of modular symbols on $\gG_0(N) \subset \SL_2(\Z)$ is defined to be the group $H^1_c(Y_0(N)(\CC),\CC)$ of compactly supported cohomology classes on the open modular curve $Y_0(N)(\CC)=\gG_0(N)\bs\cH$, where $\cH$ denotes the upper half-plane.  Let $S_2(\gG_0(N))$ denote the space of cuspidal modular forms for $\gG_0(N)$.  By the Eichler-Shimura isomorphism, the space $S_2(\gG_0(N))$ embeds into $H^1_c(Y_0(N)(\CC),\CC)$ and the image can be characterized by the action of the Hecke operators.  In sum, to compute with the space of modular forms $S_2(\gG_0(N))$, one can equivalently compute with the space of modular symbols $H^1_c(Y_0(N)(\CC),\CC)$ together with its Hecke action.  This latter space is characterized by a natural isomorphism of Hecke-modules
\[ H^1_c(Y_0(N)(\CC),\CC)\cong \Hom_{\gG_0(N)}(\Div^0\PP^1(\QQ),\CC), \]
where a cohomology class $\go$ is mapped to the linear functional which sends the divisor $s-r \in \Div^0\PP^1(\QQ)$ to the integral of $\go$ over the image on $Y_0(N)$ of a path in $\cH$ between the cusps $r$ and $s$.

Modular symbols have proved to be crucial in both computational and theoretical roles.  They arise in the study of special values of $L$-functions of classical modular forms, in the formulation of $p$-adic measures and $p$-adic $L$-functions, as well as in the conjectural constructions of Gross-Stark units~\cite{DarmonDasgupta} and Stark-Heegner points~\cite{DarmonSHP}.  It is therefore quite inconvenient that a satisfactory formalism of modular symbols is absent in the context of automorphic forms on other Shimura varieties.  Consequently, the corresponding theory is not as well understood.  From this point of view, alternative methods for the explicit study of Hilbert modular forms are of particular interest.  

Let $F$ be a totally real field of degree $n=[F:\QQ]$ and let $\ZZ_F$ denote its ring of integers.  Let $S_2(\fN)$ denote the Hecke module of (classical) Hilbert modular cusp forms over $F$ of parallel weight $2$ and level $\fN \subset \ZZ_F$.  Demb\'el\'e~\cite{Dembele} and Demb\'el\'e and Donnelly~\cite{DD} have presented methods for computing with the space $S_2(\fN)$ under the assumption that $n$ is even.  Their strategy is to apply the Jacquet-Langlands correspondence in order to identify systems of Hecke eigenvalues occurring in $S_2(\fN)$ inside spaces of automorphic forms on $B^\cp$, where $B$ is the quaternion algebra over $F$ ramified precisely at the infinite places of $F$---whence their assumption that $n$ is even.  In this way, Demb\'el\'e and his coauthors have convincingly demonstrated that automorphic forms on totally definite quaternion algebras, corresponding to Shimura varieties of dimension zero, are amenable to computation.  

Here, we provide an algorithm which permits this computation when $n$ is odd.  We locate systems of Hecke eigenvalues in the (degree one) cohomology of a Shimura curve.  We explain how a reduction theory for the associated quaternionic unit groups, arising from a presentation of a fundamental domain for the action of this group \cite{V-fd}, allows us to compute the Hecke module structure of these cohomology groups in practice.  Our methods work without reference to cusps or a canonical moduli interpretation of the Shimura curve, as these features of the classical situation are (in general) absent.

Our main result is as follows.

\begin{theorema*}
There exists an algorithm which, given a totally real field $F$ of strict class number $1$ and odd degree $n$, and an ideal $\fN$ of $\ZZ_F$, computes the system of Hecke eigenvalues associated to Hecke eigenforms in the space $S_2(\fN)$ of Hilbert modular forms of parallel weight $2$ and level $\fN$.
\end{theorema*}

In fact, our methods work more generally for fields $F$ of even degree as well (under a hypothesis on $\fN$) and for higher weight $k$, and therefore overlap with the methods of Demb\'el\'e and his coauthors in many cases; see the precise statement in Theorem \ref{maintheorem} below.  This overlap follows from the Jacquet-Langlands correspondence, but it can also be explained by the theory of nonarchimedean uniformization of Shimura curves: one can describe the $\CC_\fp$-points of a Shimura curve, for suitable primes $\fp$, as the quotient of the $\fp$-adic upper half-plane $\cH_\fp$ by a definite quaternion order of the type considered by Demb\'el\'e.  (For a discussion of the assumption on the class number of $F$, see Remark \ref{strictclassno1}.)  In particular, this theorem answers (in part) a challenge of Elkies \cite{Elkies} to compute modular forms on Shimura curves. 

The article is organized as follows.  In \S\S 1--2, we introduce Hilbert modular forms (\S\ref{S:HMF}) and quaternionic modular forms (\S\ref{S:QMF}) and the correspondence of Jacquet-Langlands which relates them.  In \S\ref{S:cohom}, we discuss how systems of Hecke eigenvalues associated to certain Hilbert or quaternionic modular forms may be found in cohomology groups of Shimura curves.  The rest of the paper is devoted to explicit computation in these cohomology groups.  In~\S\ref{S:quaternions}, we discuss algorithms for representing quaternion algebras and their orders.  In \S\ref{S:fuchsian}, we show how the fundamental domain algorithms of the second author allow for computation in the cohomology of Fuchsian groups, and we show how algorithms for solving the word problem in these groups is the key to computing the Hecke action.  We conclude by presenting applications and examples of our algorithms.

Beyond applications to the enumeration of automorphic forms, the techniques of this paper hold the promise of applications to Diophantine equations.  The first author \cite{Greenberg} has proposed a conjectural, $p$-adic construction of algebraic Stark-Heegner points on elliptic curves over totally real fields.  These points are associated to the data of an embedding of a non-CM quadratic extension $K$ of a totally real field $F$ into a quaternion $F$-algebra $B$.  Note that such a quaternion algebra cannot be totally definite.  We propose to generalize the formalism of overconvergent modular symbols employed in~\cite{DP} in the case $B=\MM_2(\QQ)$ to the general quaternionic situation in order to allow for the efficient calculation of these points.

The authors gratefully acknowledge the hospitality of the \textsf{Magma} group at the University of Sydney for their hospitality and would like to thank Lassina Demb\'el\'e, Steve Donnelly, Benjamin Linowitz, and Ron Livn\'e for their helpful comments.

\section{Hilbert modular forms}\label{S:HMF}

We begin by defining the space of classical Hilbert modular cusp forms. Our main reference for standard facts about Hilbert forms is Freitag \cite{Freitag}.

Let $F$ be a totally real field of degree $n=[F:\QQ]$ with ring of integers $\ZZ_F$.  Throughout, we assume that the strict class number of $F$ is equal to $1$ (see Remark \ref{strictclassno1} for comments on this assumption).  Let $v_1,\ldots,v_n$ be the embeddings of $F$ into $\RR$.  If $x\in F$, we write $x_i$ as a shorthand for $v_i(x)$.  Each embedding $v_i$ induces a  embedding $v_i:\MM_2(F)\hookrightarrow \MM_2(\RR)$.  Extending our shorthand to matrices, if $\gg\in \MM_2(F)$ we write $\gg_i$ for $v_i(\gg)$.  Let 
\[
\GL_2^+(F) = \{\gg\in\GL_2(F) : \text{$\det \gg_i>0$ for all $i=1,\ldots,n$}\}.
\]
The group $\GL_2^+(F)$ acts on the cartesian product $\cH^n$ by the rule
\[
\gg(\gt_1,\ldots,\gt_n)=(\gg_1\gt_1,\ldots,\gg_n\gt_n)
\]
where as usual $\GL_2^+(\R)$ acts on $\cH$ by linear fractional transformations.  

Let $\gg=\left(\begin{smallmatrix}a&b\\c&d\end{smallmatrix}\right)\in \GL_2(\RR)$ and $\gt\in\cH$.  We define
\begin{equation} \label{jlabel}
j(\gg,\gt) = c\gt+d\in\CC.
\end{equation}
For a \emph{weight}
\begin{equation}\label{E:weight}
k=(k_1,\ldots,k_n) \in (2\ZZ_{>0})^n,
\end{equation}
we define a right \emph{weight $k$ action} of $\GL_2^+(F)$ on the space of complex-valued functions on $\cH^n$ by
\[
(f\slsh{k} \gg)(\gt)=f(\gg\gt)\prod_{i=1}^n (\det\gg_i)^{k_i/2}j(\gg_i,\gt_i)^{-k_i}
\]
for $f : \cH^n \to \CC$ and $\gamma \in \GL_2^+(F)$.  The center $F^\cp$ of $\GL_2^+(F)$ acts trivially on such $f$.   Therefore, the weight $k$ action descends to an action of $\PGL_2^+(F)=\GL_2^+(F)/F^\cp$.

Now let $\fN$ be a (nonzero) ideal of $\ZZ_F$.  Define 
\[
\gG_0(\fN)=\left\{\begin{pmatrix}a&b\\c&d\end{pmatrix}\in\GL_2^+(\ZZ_F) : c\in\fN\right\}.
\]
A \emph{Hilbert modular cusp form of weight $k$ and level $\gG_0(\fN)$} is an analytic function $f:\cH^n\to\CC$ such that $f\slsh{k}\gg=f$ for all $\gg\in\gG_0(\fN)$ and such that $f$ vanishes at the cusps of $\gG_0(\fN)$.  We write $S_k(\fN)=S_k(\gG_0(\fN))$ for the finite-dimensional $\C$-vector space of Hilbert modular forms of weight $k$ and level $\gG_0(\fN)$.  (See Freitag \cite[Chapter 1]{Freitag} for proofs and a more detailed exposition.)

The space $S_k(\fN)$ is equipped with an action of Hecke operators as follows.  Let $\fp$ be a (nonzero) prime ideal of $\ZZ_F$ with $\fp\nmid \fN$.  Write $\F_\fp$ for the residue field of $\ZZ_F$ at $\fp$.  By our assumption that $F$ has strict class number one, there exists a totally positive element $p \in \ZZ_F$ which generates the ideal $\fp$.  Let $\pi=\begin{pmatrix} 1 & 0 \\ 0 & p \end{pmatrix}$.  Then there are elements $\gg_a\in \Gamma=\gG_0(\fN)$, indexed by $a \in \PP^1(\F_\fp)$, such that
\begin{equation} \label{Gammagaa}
\Gamma \pi \Gamma = \bigsqcup_{a\in\PP^1(\F_\fp)}\Gamma\ga_a,
\end{equation}
where $\ga_a = \pi \gg_a$.
If $f\in S_k(\fN)$, we define
\[
f\slsh{} T_\fp=\sum_{a\in\PP^1(\F_\fp)} f\slsh{k}\ga_a.
\]
Then $f \slsh{} T_\fp$ also belongs to $S_k(\fN)$ and the operator $T_\fp:S_k(\fN)\to S_k(\fN)$ is called the \emph{Hecke operator} associated to the ideal $\fp$.  One verifies directly that $T_\fp$ depends only on $\fp$ and not on our choice of generator $p$ of $\fp$ or on our choice of representatives (\ref{Gammagaa}) for $\Gamma\bs\Gamma\gp\Gamma$.

Suppose instead now that $\fp$ is a prime ideal of $\ZZ_F$ such that $\fp^e\parallel\fN$.  Let $p$ and $n$ be totally positive generators of $\fp$ and $\fn$, respectively.  Then there exist $x,y\in \ZZ_F$ such that $xp^e - y(n/p^e)=1$.  The element
\[
\gp=\begin{pmatrix}xp^e & y\\ n & p^e\end{pmatrix}
\]
with $\det\gp=p^e$ normalizes $\Gamma$, and we have $\gp^2\in F^\cp \Gamma$.  Setting
\[
f\slsh{} W_{\fp^e}=f\slsh{}\gp,
\]
we verify that $W_{\fp^e}$ is an involution of $S_2(\fN)$, called an \emph{Atkin-Lehner involution}, and this operator depends only on the prime power $\fp^e$ and not on its generator $p$ or the elements $x,y$.  

We conclude this section by defining the space of newforms.  Let $\fM$ be an ideal of $\Z_F$ with $\fM \mid \fN$.  For any totally positive element $d \mid \fN\fM^{-1}$ we have a map
\begin{align*}
h_d: S_k(\gG_0(\fM)) &\hookrightarrow S_k(\fN) \\
f &\mapsto f \slsh{} \left( \begin{smallmatrix} d & 0 \\ 0 & 1 \end{smallmatrix} \right).
\end{align*}
We say that $f \in S_k(\fN)$ is an \emph{oldform} at $\fM$ if $f$ is in the image of $h_d$ for some $d \mid \fN\fM^{-1}$.  Let $S_k(\fN)^{\text{$\fM$-old}}$ denote the space of oldforms at $\fM$.  Then we can orthogonally decompose the space $S_k(\fN)$ as 
\[ S_k(\fN)=S_k(\fN)^{\text{$\fM$-old}} \oplus S_k(\fN)^{\text{$\fM$-new}} \]
and we say that $f \in S_k(\fN)$ is a \emph{newform} at $\fM$ if $f \in S_k(\fN)^{\text{$\fM$-new}}$.

\section{Quaternionic modular forms}\label{S:QMF}

Our main reference for this section is Hida~\cite[\S2.3]{HidaHMF}.  Let $B$ be a quaternion algebra over $F$ which is split at the real place $v_1$ and ramified at the real places $v_2,\ldots,v_n$ (and possibly some finite places).  Let $\fD$ be the \emph{discriminant} of $B$, the product of the primes of $\Z_F$ at which $B$ is ramified.  Let $\go(\fD)$ denote the number of distinct primes dividing $\fD$.  Then since a quaternion algebra is ramified at an even number of places, we have
\begin{equation} \label{omegaD}
\go(\fD) \equiv n-1\pmod{2}.
\end{equation}
We note that the case $\fD=(1)$ is possible in (\ref{omegaD}) if and only if $n$ is odd.  For unity of presentation, we assume that $\go(\fD) + n>1$ or equivalently that $B\ncong \MM_2(\QQ)$ or that $B$ is a division algebra.

Since $B$ is split at $v_1$, we may choose an embedding
\begin{equation} \label{gidef}
\gi_1:B \hookrightarrow B\tp \RR\xrightarrow{\sim} \MM_2(\RR).
\end{equation}

We denote by $\nrd:B \to F$ the \emph{reduced norm} on $B$, defined by $\nrd(\gamma)=\gamma\overline{\gamma}$ where $\overline{\phantom{x}}:B \to B$ is the unique \emph{standard involution} (also called \emph{conjugation}) on $B$.  Let $B^\cp_+$ denote the subgroup of $B^\cp$ consisting of elements with totally positive reduced norm.  Since $B$ is ramified at all real places except $v_1$, an element $\gamma \in B^\times$ has totally positive norm if and only if $v_1(\nrd(\gamma))>0$, so that
\begin{equation} \label{Bcpplus}
B^\cp_+=\{\gg\in B^\cp : v_1(\nrd \gg) = \det \gi_1(\gg)>0\}.
\end{equation}
The group $B^\cp_+$ acts on $\cH$ via $v_1$; we write simply $\gg_1$ for $\iota_1(\gg)$.  As $F^\cp\subs B^\cp_+$ acts trivially on $\cH$ via $\gi_1$, the action of $B^\cp_+$ on $\cH$ descends to the quotient $B^\cp_+/F^\cp$.

For an integer $m\geq 0$, let $P_m=P_m(\CC)$ be the subspace of $\CC[x,y]$ consisting of homogeneous polynomials of degree $m$.  In particular, $P_0=\CC$.  For $\gg\in\GL_2(\CC)$, let $\bar{\gg}$ be the adjoint of $\gg$, so that $\gg\bar{\gg}=\det\gg$.  Note that this notation is consistent with the bar notation used for conjugation in a quaternion algebra as $\gi_1$ is an isomorphism of algebras with involution: $\gi_1(\bar{\gg})=\overline{\gi_1(\gg)}$.  Define a right action of $\GL_2(\CC)$ on $P_{m}(\CC)$ by
\[
(q \cdot \gg)(x,y) = 
q((x\,\,y)\bar{\gg})=q(dx-cy,-bx+ay)
\]
for $\gamma=\left(\begin{smallmatrix}a&b\\c&d\end{smallmatrix}\right) \in \GL_2(\CC)$ and $q \in P_m(\CC)$.  For $\ell\in\ZZ$, define a modified right action $\cdot_\ell$ of $\GL_2(\CC)$ on $P_m(\CC)$ by
\[
q\cdot_\ell \gg = (\det\gg)^{\ell}\, (q\cdot\gg).
\]
We write $P_m(\ell)(\CC)$ for the resulting right $\GL_2(\CC)$-module. 

Let $k$ be a weight as in (\ref{E:weight}) and let $w_i=k_i-2$ for $i=2,\ldots,n$.  Define the right $\GL_2(\CC)^{n-1}$-module
\begin{equation} \label{WCC}
W(\CC)= P_{w_2}(-w_2/2)(\CC)\tp\cdots\tp P_{w_n}(-w_n/2)(\CC).
\end{equation}
For the ramified real places $v_2,\dots,v_n$ of $F$, we choose splittings
\begin{equation} \label{Bcpplus2n}
\iota_i : B \hookrightarrow B \otimes_F \CC \cong M_2(\CC).
\end{equation}
We abbreviate as above $\gamma_i = \iota_i(\gamma)$ for $\gamma \in B$.  Then $W(\CC)$ becomes a right $B^\cp$-module via $\gamma \mapsto (\gamma_2,\dots,\gamma_n) \in \GL_2(\CC)^{n-1}$.  We write $(x,\gg)\mapsto x^\gg$ for this action.  We define a \emph{weight $k$ action} of $B^\cp_+$ on the space of $W(\CC)$-valued functions on $\cH$ by
\[
(f \slsh{k} \gg)(\gt)=(\det \gg_1)^{k_1/2}j(\gg_1,\gt)^{-k_1}f(\gg_1\gt)^\gg
\]
for $f:\cH \to W(\C)$ and $\gamma \in B_+^\times$, where $j$ is defined as in (\ref{jlabel}).  Note that the center $F^\cp \subset B^\cp_+$ again acts trivially, so the action descends to $B^\cp_+/F^\cp$.  We endow $W(\CC)$ with an analytic structure via a choice of linear isomorphism $W(\CC)\cong \CC^{(k_2-1)\cdots(k_n-1)}$; the resulting analytic structure does not depend on this choice.  

Let $\fN$ be an ideal of $\ZZ_F$ which is prime to $\fD$ and let $\cO_0(\fN)$ be an Eichler order in $B$ of level $\fN$.  We denote by $\cO_0(\fN)_+^\cp \subset B_+^\times$ the units of $\cO_0(\fN)$ with totally positive reduced norm, as in (\ref{Bcpplus}), and we let $\Gamma_0^{\fD}(\fN) = \cO_0(\fN)_+^\cp/\Z_{F,+}^\cp$, where $\Z_{F,+}^\cp$ denotes the units of $\Z_F$ with totally positive norm.  A \emph{quaternionic modular form of weight $k$ and level $\cO_0(\fN)$} is an analytic function $f:\cH \to W(\CC)$ such that $f \slsh{k} \gg=f$ for all $\gg\in \cO_0(\fN)^\cp_+$. 
We write $S_k^{\fD}(\fN)$ for the finite-dimensional $\C$-vector space of quaternionic modular forms of weight $k$ and level $\cO_0(\fN)$.

Spaces of quaternion modular forms can be equipped with the action of Hecke operators.  Let $\fp$ be a prime ideal of $\ZZ_F$ with $\fp\nmid \fD\fN$.  Since $F$ has strict class number $1$, by strong approximation \cite[Theor\`eme III.4.3]{Vigneras} there exists $\gp \in \cO_0(\fN)$ such that $\nrd \gp$ is a totally positive generator for the ideal $\fp$.  It follows that there are elements $\gg_a\in \Gamma = \gG_0(\fN)$, indexed by $a\in \PP^1(\F_\fp)$, such that
\begin{equation} \label{heckequat}
\Gamma\gp\Gamma = \bigsqcup_{a\in\PP^1(\F_\fp)}\Gamma\ga_a,
\end{equation}
where $\ga_a= \gp\gg_a$.  We define the \emph{Hecke operator} $T_\fp: S_k(\fN) \to S_k(\fN)$ by the rule
\begin{equation} \label{heckequatdef}
f\slsh{}T_\fp = \sum_{a\in\PP^1(\F_\fp)}f\slsh{k}\ga_a.
\end{equation}

The space $S_k^{\fD}(\fN)$ also admits an action of Atkin-Lehner operators.  Now suppose that $\fp^e \subset \ZZ_F$ is a prime power with $\fp^e \parallel \fD\fN$.  (Recall that $e=1$ if $\fp \mid \fD$ and that $\fD$ and $\fN$ are coprime.)  Then there exists an element $\pi \in \cO_0(\fN)$ whose reduced norm is a totally positive generator of $\fp^e$ and such that $\pi$ generates the unique two-sided ideal of $\cO_0(\fN)$ with norm $\fp^e$.  The element $\pi$ normalizes $\cO_0(\fN)$ and $\pi^2\in \cO_0(\fN) \subset \cO_0(\fN)^\times F^\times$ (see Vign\'eras \cite[Chapitre II, Corollaire 1.7]{Vigneras} for the case $\fp\mid\fD$ and the paragraph following \cite[Chapitre II, Lemme 2.4]{Vigneras} if $\fp\mid\fN$).  Thus, we define the \emph{Atkin-Lehner involution} $W_{\fp^e}: S_k(\fN) \to S_k(\fN)$ by
\begin{equation} \label{AtkinLehner2}
f \slsh{} W_{\fp^e}=f\slsh{k}\gp
\end{equation}
for $f \in S_k(\fN)$.
As above, this definition is independent of the choice of $\pi$ and so only depends on the ideal $\fp^e$.

The following fundamental result---the Jacquet-Langlands correspondence---relates these two spaces of Hilbert and quaternionic modular forms.

\begin{theorem} \label{JLthm}
There is a vector space isomorphism
\[
S_k^\fD(\fN) \xrightarrow{\sim} S_k(\fD\fN)^{\text{$\fD$-new}}
\]
which is equivariant for the actions of the Hecke operators $T_\fp$ with $\fp\nmid\fD\fN$ and the Atkin-Lehner involutions $W_{\fp^e}$ with $\fp^e\parallel\fD\fN$.
\end{theorem}

A useful reference for the Jacquet-Langlands correspondence is Hida~\cite[Proposition 2.12]{HidaCM}, where Theorem~\ref{JLthm} is deduced from the representation theoretic results of Jacquet-Langlands.  

In particular, when $n=[F:\Q]$ is odd, we may take $\fD=(1)$ to obtain an isomorphism $S_k^{(1)}(\fN) \xrightarrow{\sim} S_k(\fN)$.  

\section{Quaternionic modular forms and the cohomology of Shimura curves}\label{S:cohom}

In this section, we relate the spaces $S_k^{\fD}(\fN)$ of quaternionic modular forms, together with their Hecke action, to the cohomology of Shimura curves.  As above, let $\gG=\gG_0^{\fD}(\fN)=\cO_0(\fN)^\cp_+/\ZZ_F^\cp$.  The action of $\gG$ on $\cH$ is properly discontinuous.  Therefore, the quotient
$\Gamma\bs\cH$ has a unique complex structure such that the natural projection $\cH\to \gG\bs\cH$ is analytic.  Since $B$ is, by assumption, a division algebra, $\gG\bs\cH$ has the structure of a compact Riemann surface.  By the theory of canonical models due to Shimura \cite{ShimuraZeta} and Deligne \cite{Deligne}, this Riemann surface is the locus of complex points of the \emph{Shimura curve} $X=X_0^{\fD}(\fN)$ which is defined over $F$, under our assumption that $F$ has strict class number $1$.

Define the right $\GL_2(\CC)^n=\GL_2(\CC)\cp\GL_2(\CC)^{n-1}$-module
\[
V(\CC)=\bigotimes_{i=1}^{n} P_{w_i}(-w_i/2)(\CC) =P_{w_1}(-w_1/2)(\CC)\tp W(\CC).
\]
The group $B^\cp$ acts on $V(\CC)$ via the embedding $B^\cp\hookrightarrow\GL_2(\CC)^n$ given by $\gg\mapsto(\gg_1,\ldots\gg_n)$ arising from (\ref{gidef}) and (\ref{Bcpplus2n}).
As $F^\cp$ acts trivially on $P_{w_1}(-w_1/2)(\CC)$ via $\gi_1$ and on $W(\CC)$ via $(\gi_2,\ldots,\gi_n)$, the action of $\cO_0(\fN)^\cp_+ \subset B^\cp$ on $V(\CC)$ descends to a right action of $\gG$ on $V(\CC)$.

It is convenient to identify $V(\CC)$ with the subspace of the algebra $\CC[x_1,y_1,\ldots,x_n,y_n]$ consisting of those polynomials $q$ which are homogeneous in $(x_i,y_i)$ of degree $w_i$.  Under this identification, the action of $B^\cp$ on $V(\CC)$ takes the form
\[
q^\gg(x_1,y_1,\ldots,x_n,y_n)= \prod_{i=1}^n(\det\gg_i)^{-w_i/2}
q((x_1\,\,y_1)\bar{\gg}_1,\ldots(x_n\,\,y_n)\bar{\gg}_n).
\]

Let $V=V(\CC)$.  We now consider the cohomology group $H^1(\gG,V)$.  We represent elements of $H^1(\Gamma,V)$ by (equivalence classes of) crossed homomorphisms.  Recall that a \emph{crossed homomorphism} (or \emph{$1$-cocycle}) $f:\Gamma\to V$ is a map which satisfies the property 
\begin{equation} \label{crossedprop}
f(\gg\delta)=f(\gg)^\delta + f(\delta)
\end{equation}
for all $\gamma,\delta \in \Gamma$.  A crossed homomorphism $f$ is \emph{principal} (or a \emph{$1$-coboundary}) if there is an element $g\in V$ such that 
\begin{equation} \label{principalprop}
f(\gg) = g^{\gg}- g
\end{equation} for all $\gamma \in \Gamma$.  Let $Z^1(\Gamma,V)$ and $B^1(\Gamma,V)$ denote the spaces of crossed homomorphisms and principal crossed homomorphisms from $\Gamma$ into $V$, respectively.  Then
\[
H^1(\Gamma,V)=Z^1(\Gamma,V)/B^1(\Gamma,V).
\]

We now define (for a third time) the action of the Hecke operators, this time in cohomology.  Let $f:\Gamma\to V$ be a crossed homomorphism, and let $\gg\in \Gamma$.   Recall the definition of Hecke operators (\ref{heckequat}--\ref{heckequatdef}).
There are elements $\delta_a \in \Gamma$ for $a \in \PP^1(\F_\fp)$ and a unique permutation $\gg^*$ of $\PP^1(\F_\fp)$ such that
\begin{equation} \label{deltaa}
\ga_a\gg = \delta_a\ga_{\gg^*a}
\end{equation}
for all $a$.  Define $f \slsh{} T_\fp:\Gamma\to V$ by
\begin{equation}\label{E:heckeformula}
(f \slsh{} T_\fp)(\gg) = \sum_{a\in\PP^1(\F_\fp)}f(\delta_a)^{\ga_a}.
\end{equation}
It is a standard calculation~\cite[\S8.3]{ShimuraArith} that $f \slsh{} T_\fp$ is a crossed homomorphism and that $T_\fp$ preserves $1$-coboundaries.  Moreover, $f\slsh{} T_\fp$ does not depend on our choice of the coset representatives $\ga_a$.  Therefore, $T_\fp$ yields a well-defined operator
\[
T_\fp:H^1(\Gamma,V)\to H^1(\Gamma,V).
\]

\begin{remark}
One may in a natural way extend this definition to compute the action of any Hecke operator $T_\fa$ for $\fa\subset \Z_F$ an ideal with $(\fa, \fD\fN)=1$.  
\end{remark}

We now define an operator corresponding to complex conjugation.  Let $\mu \in \cO_0(\fN)^\cp$ be an element such that $v(\nrd \mu)<0$; such an element exists again by strong approximation.  Then $\mu$ normalizes $\Gamma$ and $\mu^2 \in \gG$.  If $f\in Z^1(\Gamma,V)$ then the map $f \slsh{} W_{\infty}$ defined by
\begin{equation} \label{winfty}
(f \slsh{} W_{\infty})(\gg) = f(\mu \gamma \mu^{-1})^{\gm}
\end{equation}
is also a crossed homomorphism which is principal if $f$ is principal, so it induces a linear operator
\[
W_\infty:H^1(\Gamma,V)\to H^1(\Gamma,V).
\]
Since $\mu^2\in \gG$ and $\Z_{F,+}^*$ acts trivially on $V$, the endomorphism $W_\infty$ has order two.  Therefore, $H^1(\Gamma,V)$ decomposes into eigenspaces for $W_{\infty}$ with eigenvalues $+1$ and $-1$, which we denote
\[
H^1(\Gamma,V) = H^1(\Gamma,V)^+\ds H^1(\Gamma,V)^-.
\]
It is not hard to see that $T_\fp$ commutes with $W_{\infty}$.  Therefore, $T_\fp$ perserves the $\pm$-eigenspaces of $H^1(\Gamma,V)$.

The group $H^1(\Gamma,V)$ also admits an action of Atkin-Lehner involutions.
Letting $\fp$, $e$, and $\gp$ be as in the definition of the Atkin-Lehner involutions (\ref{AtkinLehner2}) in \S\ref{S:QMF}, we define the involution $W_{\fp^e}$ on $H^1(\Gamma, V)$ by
\begin{equation} \label{AL3}
(f \slsh{} W_{\fp^e})(\gamma)=f(\pi \gamma \pi^{-1})^{\gp}.
\end{equation}

We conclude this section by relating the space of quaternionic modular forms to cohomology by use of the \emph{Eichler-Shimura isomorphism} (in analogy with the classical case $B=\M_2(\QQ)$ \cite[\S8.2]{ShimuraArith}).  We choose a base point $\gt\in\cH$, and for $f\in S_k^{\fD}(\fN)$ we define a map
\[
\ES(f):\Gamma\to V
\]
by the rule
\[
\ES(f)(\gg) = \int_{\gg_1^{-1}\gt}^\gt f(z)(zx_1+y_1)^{w_1}dz\in V(\CC).
\]
A standard calculation~\cite[\S8.2]{ShimuraArith} shows that $\ES(f)$ is a crossed homomorphism which depends on the choice of base point $\gt$ only up to $1$-coboundaries.  Therefore, $\ES$ descends to a homomorphism
\[
\ES:S_k(\Gamma)\to H^1(\Gamma,V).
\]
Let
\[
\ES^\pm:S_k(\Gamma)\to H^1(\Gamma,V)^\pm.
\]
be the composition of $\ES$ with projection to the $\pm$-subspace, respectively.

\begin{theorem}[{\cite[\S4]{MS}}] \label{EichlerShimurathm}
The maps $\ES^\pm$ are isomorphisms (of $\C$-vector spaces) which are equivariant with respect to the action of the Hecke operators $T_\fp$ for $\fp \nmid \fD\fN$ and $W_{\fp^e}$ for $\fp^e\parallel\fD\fN$.
\end{theorem}

In sum, to compute the systems of Hecke eigenvalues occurring in spaces of Hilbert modular forms, combining the Jacquet-Langlands correspondence (Theorem \ref{JLthm}) and the Eichler-Shimura isomorphism (Theorem \ref{EichlerShimurathm}), it suffices to enumerate those systems occurring in the Hecke module $H^1(\Gamma,V(\CC))^+$.

Let $\fN\subs\ZZ_F$ be an ideal.  Suppose that $\fN$ admits a factorization $\fN=\fD\fM$ such that
\begin{equation} \label{factorz}
\text{$\fD$ is squarefree, $\fD$ is coprime to $\fM$, and $\go(\fD)\equiv n-1\pmod{2}$.}
\end{equation}
Then there is a quaternion $F$-algebra $B$ ramified precisely at primes dividing $\fD$ and at the infinite places $v_2,\ldots,v_n$ of $F$.  The goal of the second part of this paper is to prove the following theorem which describes how spaces of automorphic forms discussed in the first part of this paper can be computed in practice.

\begin{theorem} \label{maintheorem}
There exists an explicit algorithm which, given a totally real field $F$ of strict class number $1$, an ideal $\fN\subs\ZZ_F$, a factorization $\fN=\fD\fM$ as in \textup{(\ref{factorz})}, and a weight $k \in (2\Z_{>0})^n$, computes the systems of eigenvalues for the Hecke operators $T_\fp$ with $\fp\nmid \fD\fM$ and the Atkin-Lehner involutions $W_{\fp^e}$ with $\fp^e\parallel\fD\fM$ which occur in the space of Hilbert modular cusp forms $S_k(\fN)^{\text{$\fD$-new}}$.
\end{theorem}

By this we mean that we will exhibit an explicit finite procedure which takes as input the field $F$, the ideals $\fD$ and $\fM$, and the integer $k$, and produces as output a set of sequences $(a_\fp)_{\fp}$ with $a_\fp \in \overline{\QQ}$ which are the Hecke eigenvalues of the finite set of $\fD$-newforms in $S_2(\fN)$.  The algorithm will produce the eigenvalues $a_\fp$ in any desired ordering of the primes $\fp$.  

\begin{remark} \label{strictclassno1}
Generalizations of these techniques will apply when the strict class number of $F$ is greater than $1$.  In this case, the canonical model $X$ for the Shimura curve associated to $B$ is the disjoint union of components defined over the strict class field of $F$ and indexed by the strict ideal class group of $\Z_F$.  The Hecke operators and Atkin-Lehner involutions then permute these components nontrivially and one must take account of this additional combinatorial data when doing computations.  Because of this additional difficulty (see also Remark \ref{oopsfunddomtoo} below), we leave this natural extension as a future project.  

We note however that it is a folklore conjecture that if one orders totally real fields by their discriminant, then a (substantial) positive proportion of fields will have strict class number $1$.  For this reason, we are content to considering a situation which is already quite common.
\end{remark}

\section{Algorithms for quaternion algebras}\label{S:quaternions}

We refer to work of Kirschmer and the second author \cite{KV} as a reference for algorithms for quaternion algebras.  We will follow the notation and conventions therein.  

To even begin to work with the algorithm implied by Theorem \ref{maintheorem} we must first find a representative quaternion algebra $B=\quat{a,b}{F}$ with discriminant $\fD$ which is ramified at all but one real place.  From the point of view of effective computability, one can simply enumerate elements $a,b \in \ZZ_F \setminus \{0\}$ and then compute the discriminant of the corresponding algebra \cite{VoightMatrixRing} until an appropriate representative is found.  (Since such an algebra exists, this algorithm always terminates after a finite amount of time.)  In practice, it is much more efficient to compute as follows.  

\begin{alg}
This algorithm takes as input a discriminant ideal $\fD \subset \Z_F$, the product of distinct prime ideals with $\omega(\fD) \equiv n-1 \pmod{2}$, and returns $a,b \in \Z_F$ such that the quaternion algebra $B=\quat{a,b}{F}$ has discriminant $\fD$.

\begin{enumalg}
\item Find $a \in \fD$ such that $v(a)>0$ for at most one real place $v_1$ of $F$ and such that $a\ZZ_F=\fD\fb$ with $\fD+\fb=\ZZ_F$ and $\fb$ odd.  
\item Find $t \in \ZZ_F/8a\ZZ_F$ such that the following hold:
\begin{itemize}
\item For all primes $\fp \mid \fD$, we have $\displaystyle{\legen{t}{\fp}=-1}$;
\item For all primes $\fq \mid \fb$, we have $\displaystyle{\legen{t}{\fq}=1}$; and
\item For all prime powers $\fr^e \parallel 8\ZZ_F$ with $\fr \nmid \fD$, we have $t \equiv 1 \pmod{\fr^e}$.
\end{itemize}
\item Find $m \in \ZZ_F$ such that $b=t+8am \in \ZZ_F$ is prime and such that $v(b)<0$ for all $v \neq v_1$ and either $v_1(a)>0$ or $v_1(b)>0$.  
\end{enumalg}
\end{alg}

Since our theorem does not depend on the correctness of this algorithm, we leave it to the reader to verify that the algebra $B=\quat{a,b}{F}$ output by this algorithm has indeed the correct set of ramified places.

\begin{remark}
When possible, it is often helpful in practice to have $a\ZZ_F=\fD$, though this requires the computation of a generator for the ideal $\fD$ with the correct real signs; for example, if $\fD=\ZZ_F$ and there exists a unit $u \in \ZZ_F^\cp$ such that $v(u)>0$ for a unique real place $v=v_1$ and such that $u \equiv 1 \pmod{8}$, then we may simply take $B=\quat{-1,u}{F}$.  

One may wish to be alternate between Steps 2 and 3 in searching for $b$.  Finally, we note that in Step 2 one may either find the element $t$ by deterministic or probabilistic means.
\end{remark}

Given this representative algebra $B$, there are algorithms \cite{VoightMatrixRing} to compute a maximal order $\cO \subset B$, which is represented in bits by a pseudobasis.  Furthermore, given a prime $\fp \nmid \fD$ there exists an algorithm to compute an embedding $\iota_\fp:\cO \hookrightarrow \MM_2(\ZZ_{F,\fp})$ where $\ZZ_{F,\fp}$ denotes the completion of $F$ at $\fp$.  As a consequence, one can compute an Eichler order $\cO_0(\fN)\subs\cO$ for any level $\fN \subset \Z_F$.

\section{Computing in the cohomology of a Shimura curve} \label{S:fuchsian}

In this section, we show how to compute explicitly in first cohomology group of a Shimura curve, equipped with its Hecke module structure.  Throughout, we abbreviate $\Gamma=\Gamma_0^{\fD}(\fN)$ and $\cO=\cO_0(\fN)$ as above.

In the choices of embeddings $\iota_1,\dots,\iota_n$ in (\ref{Bcpplus}) and (\ref{Bcpplus2n}), we may assume without loss of generality that their image is contained in $\MM_2(K)$ with $K \hookrightarrow \CC$ a number field containing $F$; in other words, we may take $K$ to be any Galois number field containing $F$ which splits $B$.  In particular, we note that then that image of $\iota_1(B)$ is contained in $K \cap \RR$.  So throughout, we may work throughout with the coefficient module 
\[
V(K)=\TP_{i=1}^n P_{w_i}(K)(-w_i/2)
\] 
since obviously $V(K) \otimes_K \CC = V(\CC)$.  The $K$-vector space $V(K)$ is then equipped with an action of $\Gamma$ via $\iota$ which can be represented using exact arithmetic over $K$.  

\begin{remark}
In fact, one can take the image of $\gi$ to be contained in $\M_2(\ZZ_K)$ for an appropriate choice of $K$, where $\ZZ_K$ denotes the ring of integers of $K$.  Therefore, one could alternatively work with $V(\ZZ_K)$ as a $\ZZ_K$-module and thereby recover from this integral stucture more information about the structure of this module.  In the case where $k=(2,2,\ldots,2)$, we may work already with the coefficient module $\ZZ$, and we do so in Algorithm \ref{weight2} below.
\end{remark}

The first main ingredient we will use is an algorithm of the second author \cite{V-fd}.  

\begin{prop} \label{funddom}
There exists an algorithm which, given $\Gamma$, returns a finite presentation for $\gG$ with a minimal set of generators and a solution to the word problem for the computed presentation.
\end{prop}

\begin{remark}\label{R:sidepairing}
One must choose a point $p \in \cH$ with trivial stabilizer $\gG_p=\{1\}$ in the above algorithm, but a random choice of a point in any compact domain will have trivial stabilizer with probability $1$.  The algorithm also yields a fundmental domain $\cF$ for $\gG$ which is, in fact a hyperbolic polygon.  For each element $\gg$ of our minimal set of generators, $\cF$ and $\gg\cF$ share a single side.  In other words, there is a unique side $s$ of $\cF$ such $\gg s$ is also a side of $\cF$.  We say that $\gg$ is a \emph{side pairing element} for $s$ and $\gg s$.
\end{remark}

\begin{remark} \label{oopsfunddomtoo}
The algorithms of the second author~\cite{V-fd} concern the computation of a fundamental domain for $\cO^\cp_1/\{\pm 1\}$ acting on $\cH$, where $\cO^\cp_1$ denote the subgroup of $\cO^\cp$ consisting of elements of reduced norm $1$.  Since we have assumed that $F$ has narrow class number $1$, the natural inclusion $\cO^\cp_1 \hookrightarrow \cO^\cp_+$ induces an isomorphism
\[
\cO^\cp_1/\{\pm 1\}\xrightarrow{\sim} \cO^\cp_+/\ZZ_F^\cp=\gG.
\]
Indeed, if $\gg\in\cO_0(\fN)^\cp_+$, then $\nrd\gg \in \Z_{F,+}^\cp=\ZZ_F^{\cp 2}$ so if $\nrd\gg=u^2$ then $\gg/u$ maps to $\gg$ in the above map.
\end{remark}

Let $G$ denote a set of generators for $\gG$ and $R$ a set of relations in the generators $G$, computed as in Proposition \ref{funddom}.  We identify $Z^1(\Gamma,V(K))$ with its image under the inclusion
\begin{align*} 
j_G:Z^1(\Gamma,V(K)) &\to \bigoplus_{g \in G} V(K) \\
f &\mapsto (f(g))_{g\in G}.
\end{align*}
It follows that $Z^1(\Gamma,V(K))$ consists of those $f \in \bigoplus_{g \in G} V(K)$ which satisfy $f(r)=0$ for $r \in R$; these become linear relations written out using the crossed homomorphism property (\ref{crossedprop}), and so an explicit $K$-basis for $Z^1(\Gamma,V(K))$ can be computed using linear algebra.  The space of principal crossed homomorphisms (\ref{principalprop})
is obtained similarly, where a basis is obtained from any basis for $V(K)$.   We obtain from this a $K$-basis for the quotient
\[ 
H^1(\Gamma,V(K)) = Z^1(\Gamma,V(K))/B^1(\Gamma,V(K))
\]
and an explicit $K$-linear map $Z^1(\Gamma, V(K)) \to H^1(\Gamma,V(K))$.

We first decompose the space $H^1(\Gamma,V(K))$ into $\pm$-eigenspaces for complex conjugation $W_\infty$ as in (\ref{winfty}).  An element $\mu \in \cO^\times$ with $v_1(\nrd(\mu))<0$ can be found simply by enumeration of elements in $\calO$.  Given such an element $\mu$, we can then find a $K$-basis for the subspace $H^1(\Gamma,V(K))^+$ by linear algebra.

\begin{remark}
This exhaustive search in practice benefits substantially from the methods of Kirschmer and the second author \cite{KV} using the absolute reduced norm on $\calO$, which gives the structure of a lattice on $\calO$ so that an LLL-lattice reduction can be performed.  
\end{remark}

Next, we compute explicitly the action of the Hecke operators on the $K$-vector space $H^1(\Gamma,V)^+$.  Let $\fp\subs\ZZ_F$ be an ideal with $\fp\nmid\fD\fN$ and let $\gp\subs\cO$ be such that $\nrd\gp$ is a totally positive generator of $\gp$.  We need to compute explicitly a coset decomposition as in (\ref{heckequat}).
Now the set $\cO(\fp)=\cO^\times \pi \cO^\times$ is in natural bijection with the set of elements whose reduced norm generates $\fp$; associating to such an element the left ideal that it generates gives a bijection between the set $\cO^\times \backslash \cO(\fp)$ and the set of left ideals of reduced norm $\fp$, and in particular shows that the decomposition~\eqref{heckequat} is independent of $\pi$.  But this set of left ideals in turn is in bijection \cite[Lemma 6.2]{KV} with the set $\PP^1(\F_\fp)$: explicitly, given a splitting $\gi_\fp:\cO \hookrightarrow \MM_2(\ZZ_{F,\fp})$, the left ideal corresponding to a point $a=(x:y)\in\PP^1(\F_\fp)$ is
\begin{equation} \label{Iadef}
I_a:=\cO\gi_\fp^{-1}\begin{pmatrix} x & y \\ 0 & 0 \end{pmatrix} + \cO\fp.
\end{equation}
For $a \in \PP^1(\F_\fp)$, we let $\alpha_a \in \calO$ be such that $\calO \alpha_a = I_a$ and $\nrd \alpha_a = \nrd \pi$.  We have shown that
\[ \cO(\fp)=\cO^\times \pi \cO^\times = \bigsqcup_{a \in \PP^1(\F_\fp)} \cO^\times \alpha_a. \]
To compute the $\ga_a$, we use the following proposition {\cite{KV}}.

\begin{prop} \label{KVprinc}
There exists an (explicit) algorithm which, given a left $\cO$-ideal $I$, returns $\alpha \in \cO$ such that $\cO \alpha = I$.
\end{prop}
If $v_1(\nrd\ga_a)$ happens to be negative, then replace $\ga_a$ by $\gm\ga_a$.  
Intersection with $\cO^\cp_+$ then yields the decomposition~\eqref{heckequat}.

Next, in order to use equation~\eqref{E:heckeformula} to compute $(f \slsh{} T_\fp)(\gg)$ for all $\gg\in G$, we need to compute to compute the permutations $a\mapsto\gg^*a$ of $\PP^1(\F_\fp)=\F_\fp \cup \{\infty\}$.

\begin{alg} \label{ba}
Given $\gg\in\gG$, this algorithm returns the permutation $a\mapsto \gg^*a$.
\begin{enumalg}
\item Let $\beta=\gi_\fp(\ga_a\gamma) \in \MM_2(\ZZ_{F,\fp})$ and let $\beta_{ij}$ denote the $ij$th entry of $\beta$ for $i,j=1,2$.
\item If $\ord_\fp(\beta_{11}) \leq 0$ then return $\beta_{12}/\beta_{11} \bmod{\fp}$.
\item If $\ord_\fp(\beta_{12})=0$ or $\ord_{\fp}(\beta_{21}) > 0$ then return $\infty$.
\item Otherwise, return $\beta_{22}/\beta_{21} \bmod{\fp}$.
\end{enumalg}
\end{alg}

The proof that this algorithm gives correct output is straightforward.  

Having computed the permutation $\gg^*$, for each $a \in \PP^1(\FF_p)$ we compute $\delta_a=\ga_a\gg\ga_{\gg^*a}^{-1} \in \Gamma$ as in (\ref{deltaa}).  By Proposition \ref{funddom}, we can write $\delta_a$ as a word in the generators $G$ for $\gG$, and using the crossed homomorphism property (\ref{crossedprop}) for each $f \in Z^1(\Gamma,V(K))^+$ we compute $f \slsh{} T_\fp \in Z^1(\Gamma,V(K))^+$ by computing $(f \slsh{} T_\fp)(\gg) \in V(K)$ for $\gg\in G$.

Finally, we compute the Atkin-Lehner involutions $W_{\fp^e}$ for $\fp^e \parallel \fD\fN$ as in (\ref{AL3}).  We compute an element $\pi \in \calO$ with totally positive reduced norm such that $\pi$ generates the unique two-sided ideal $I$ of $\calO_0(\fN)$ of norm $\fp^e$.  The ideal $I$ can be computed easily \cite{KV}, and a generator $\pi$ can be computed again using Proposition \ref{KVprinc}.

We then decompose the space $H^1(\Gamma,V(K))$ under the action of the Hecke operators into Hecke irreducible subspaces for each operator $T_\fp$, and from this we compute the systems of Hecke eigenvalues using straightforward linear algebra over $K$.  (Very often it turns out in practice that a single operator $T_\fp$ is enough to break up the space into Hecke irreducible subspaces.)

For concreteness, we summarize the algorithm in the simplest case of parallel weight $k=(2,2,\ldots,2)$.  Here, we may work simply the coefficient module $\Z$ since the $\Gamma$-action is trivial.  Moreover, the output of the computation of a fundamental domain provided by Proposition \ref{funddom} is a minimal set of generators and relations with the following properties \cite[\S 5]{V-fd}: each generator $g \in G$ is labelled either elliptic or hyperbolic, and each relation becomes trivial in the free abelian quotient.

\begin{alg} \label{weight2}
Let $\fp \subset \Z_F$ be coprime to $\fD$.  This algorithm computes the matrix of the Hecke operator $T_\fp$ acting on $H^1(\Gamma,\Z)^+$.
\begin{enumalg}
\item Compute a minimal set of generators for $\Gamma$ using Proposition \ref{funddom}, and let $G$ denote the set of nonelliptic generators.  Let $H = \bigoplus_{\gamma \in G} \Z \gamma$.
\item Compute an element $\mu$ as in (\ref{winfty}) and decompose $H$ into $\pm$-eigenspaces $H^{\pm}$ for $W_\infty$.
\item Compute $\alpha_a$ for $a \in \PP^1(\F_\fp)$ as in (\ref{deltaa}) by Proposition \ref{KVprinc}.
\item For each $\gamma \in G$, compute the permutation $\gamma^*$ of $\PP^1(\F_\fp)$ using Algorithm \ref{ba}.
\item Initialize $T$ to be the zero matrix acting on $H$, with rows and columns indexed by $G$.
\item For each $\gamma \in G$ and each $a \in \PP^1(\F_\fp)$, let $\delta_a := \alpha_a \gamma \alpha_{\gamma^*(a)}^{-1} \in \Gamma$.  Write $\delta_a$ as a word in $G$ using Proposition \ref{funddom}, and add to the column indexed by $\gamma$ the image of $\delta_a \in H$.
\item Compute the action $T^+$ of the matrix $T$ on $H^+$ and return $T^+$.
\end{enumalg}
\end{alg}

We note that Step 1 is performed as a precomputation step as it does not depend on the ideal $\fp$.

\section{Higher level and relation to homology}

In this section, we consider two topics which may be skipped on a first reading.  The first considers a simplification if one fixes the quaternion algebra and varies the level $\fN$.  The second considers the relationship between our cohomological method and the well-known method of modular symbols.

\subsection*{Higher level}
Let $\cO(1)$ be a maximal order of $B$ which contains $\cO_0(\fN)$ and let $\gG(1)=\cO(1)^\cp_+/\ZZ_F^\cp$.
From Shapiro's lemma, we have an isomorphism
\[ H^1(\Gamma(1),\Coind_{\Gamma}^{\Gamma(1)} V(K) \cong H^1(\Gamma, V(K)) \]
for any $K[\Gamma]$-module $V$ where $\Coind_{\Gamma}^{\Gamma(1)} V$ denotes the coinduced module $V$ from $\Gamma$ to $\Gamma(1)$.  In particular, if one fixes $\fD$ and wishes to compute systems of Hecke eigenvalues for varying level $\fN$, one can vary the coefficient module instead of the group and so the fundamental domain algorithm need only be called once to compute a presentation for $\Gamma(1)$.

To compute the action of $\Gamma(1)$ on the coinduced module then, we need first to enumerate the set of cosets $\Gamma \backslash \Gamma(1)$.  We use a set of side pairing elements \cite{V-fd} for $\Gamma(1)$, which are computed as part of the algorithm to compute a presentation for $\Gamma$.  (Side pairing elements are defined in Remark~\ref{R:sidepairing}.)

\begin{alg} \label{cosets}
Let $\fN$ be a level and let $G$ be a set of side pairing elements for $\Gamma(1)$.  This algorithm computes representatives $\alpha \in \Gamma(1)$ such that $\Gamma \backslash \Gamma(1) = \bigsqcup \Gamma\alpha$.
\begin{enumalg}
\item Initialize $\calF := \{1\}$ and $A := \{\}$.  Let $G^{\pm} := G \cup G^{-1}$.
\item Let 
\[ \calF := \{g\gamma : g \in G,\ \gamma \in \calF,\ g\gamma\alpha^{-1} \not\in \Gamma \text{ for all $\alpha \in A$}\} \]
and let $A := A \cup \calF$.  If $\calF = \emptyset$, then return $A$; else, return to Step 2.
\end{enumalg} 
\end{alg}

\begin{proof}[Proof of correctness]
Consider the (left) Cayley graph of $\Gamma$ on the set $G^{\pm}$: this is the graph with vertices indexed by the elements of $\Gamma$ and directed edges $\gamma \to \delta$ if $\delta=g\gamma$ for some $g \in G^{\pm}$.  Then the set $\Gamma \backslash \Gamma(1)$ is in bijection with the set of vertices of this graph under the relation which identifies two vertices if they are in the same coset modulo $\Gamma$.  Since the set $G^{\pm}$ generates $\Gamma$, this graph is connected, and Step 2 is simply an algorithmic way to traverse the finite quotient graph.
\end{proof}

\begin{remark}
Although it is tempting to try to obtain a set of generators for $\Gamma$ from the above algorithm, we note that a presentation for $\Gamma$ may be arbitrarily more complicated than that of $\Gamma(1)$ and require many more elements than the number of cosets (due to the presence of elliptic elements).  For this reason, it is more efficient to work with the induced module than to compute separately a fundamental domain for $\Gamma$.
\end{remark}

To compute the Hecke operators as in \S 5 using this simplification, the same analysis applies and the only modification required is to ensure that the elements $\alpha_a$ arising in (\ref{deltaa}) satisfy $\lambda_a \in \calO_0(\fN)$---this can obtained by simply multiplying by the appropriate element $\alpha \in \Gamma \backslash \Gamma(1)$, and in order to do this efficiently one can use a simple table lookup.

\subsection*{Homology}

We now relate the cohomological approach taken above to the method using homology.  Although the relationship between these two approaches is intuitively clear as the two spaces are dual, this alternative perspective also provides a link to the theory of modular symbols, which we briefly review now.  Let $S_2(N)$ denote the $\CC$-vector space of classical modular cusp forms of level $N$ and weight $2$.    Integration defines a nondegenerate Hecke-equivariant pairing between $S_2(N)$ and the homology group $H_1(X_0(N),\ZZ)$ of the modular curve $X_0(N)$.    Let $\cS_2(N)$ denote the space of cuspidal modular symbols, i.e.\ linear combination of paths in the completed upper half-plane $\cH^*$ whose endpoints are cusps and whose images in $X_0(N)$ are a linear combination of loops.  Manin showed that there is a canonical isomorphism
\begin{center}
$\cS_2(N) \cong H_1(X_0(N),\ZZ)$.
\end{center}
If $SL_2(\ZZ)=\bigsqcup_i \Gamma_0(N) \gamma_i$,  then the set of symbols $\gamma_i\{0,\infty\}=\{\gamma_i(0),\gamma_i(\infty)\}$ generate the space $\cS_2(N)$.  We have an explicit description of the action of the Hecke operators on the space $\cS_2(N)$, and the \emph{Manin trick} provides an algorithm for writing an arbitrary modular symbol as a $\ZZ$-linear combination of the symbols $\gamma_i\{0,\infty\}$.  

The Shimura curves $X=X_0^{\fD}(\fN)$ do not have cusps, and so this method does not generalize directly.  Consider instead the sides of a 
 fundamental domain $D$ for $\Gamma=\Gamma_0^{\fD}(\fN)$, and let $G$ be the corresponding set of side pairing elements and $R$ the set of minimal relations among them \cite{V-fd}.  The side pairing gives an explicit characterization of the gluing relations which describe $X$ as a Riemann surface, hence one obtains a complete description for the homology group $H_1(X,\ZZ)$.  Let $V$ be the set of midpoints of sides of $D$.  Then for each $v \in V$, there is a unique $\gamma \in G$ such that $w=\gamma v \in V$, and the path from $v$ to $\gamma v$ in $\cH$, which we denote by $p(\gamma)=\{v,\gamma v\}$, projects to a loop on $X$.  Each relation $r=\gamma_1 \gamma_2 \cdots \gamma_t=1$ from $R$ induces the relation in the homology group $H_1(X,\ZZ)$
\begin{align} \label{indrelat}
0 &= p(1)=p(\gamma_1 \cdots \gamma_t) \notag \\
&=\{v,\gamma_1 v\} + \{\gamma_1 v,\gamma_1\gamma_2 v\} + \dots + \{\gamma_1 \cdots \gamma_{t-1} v,\gamma_1 \cdots \gamma_{t-1}\gamma_t v\} \\
&=\{v,\gamma_1 v\} + \{v,\gamma_2 v\} + \dots + \{v,\gamma_t v\} = p(\gamma_1) + p(\gamma_2) + \dots + p(\gamma_t) \notag
\end{align}
In particular, if $\gamma \in G$ is an elliptic element then $p(\gamma)=0$ in $H^1(X,\QQ)$.  Let $\cS_2^{\fD}(\fN)$ be the $\QQ$-vector space generated by $p(\gamma)$ for $\gamma \in G$ modulo the relations (\ref{indrelat}) with $r \in R$.  It follows that
\[ H_1(X,\QQ) \cong H_1(\Gamma,\QQ) \cong \cS_2^{\fD}(\fN) \]
and we call $\cS_2^{\fD}(\fN)$ the space of \emph{Dirichlet-modular symbols} for $X$ (relative to $D$).  

The Hecke operators act in an analogous way, as follows.  If $\alpha_a$ for $a \in \PP^1(\F_\fp)$ is a set of representatives as in (\ref{deltaa}), then $\alpha_a$ acts on the path $p(\gamma)=\{v,\gamma v\}$ by $\alpha_a\{v,\gamma v\}=\{\alpha_a v, \alpha_a\gamma v\}$; if $\alpha_a \gamma = \delta_a \alpha_{\gamma^* a}= \delta_a \alpha_b$ as before, then in homology we obtain
\begin{align*}
\sum_{a \in \PP^1(\F_\fp)} \{\alpha_a v, \alpha_a \gamma v\} &=
\sum_a \{\alpha_a v, \delta_a \alpha_b v\} = \sum_a \left(\{\alpha_a v, v\} + \{v,\delta_a v\} + \{\delta_a v, \delta_a \alpha_b v\}\right) \\
&= \sum_a \{v,\delta_a v\} + \sum_a \left(\{\alpha_a v, v\} + \delta_a\{v, \alpha_b v\}\right)  \\
&= \sum_a \{v,\delta_a v\} + \sum_a \left(- \{v,\alpha_a v\} + \{v, \alpha_b v\}\right) = \sum_a \{v,\delta_a v\}.
\end{align*}
Thus, the action of the Hecke operators indeed agrees with that in cohomology, and so one could also rephrase our methods in terms of Dirichlet-modular symbols.

The analogue of the Manin trick in our context is played by the solution to the word problem in $\Gamma$.  When $\Gamma=\SL_2(\Z)$, the Manin trick arises directly from the Euclidean algorithm; therefore, our methods may be seen in this light as a generalization of the Euclidean algorithm to the group $\Gamma$.  Already this point of view has been taken by Cremona \cite{CremonaIQF} and his students, who generalized methods of modular symbols to the case of $SL_2(\Z_K)$ for $K$ an imaginary quadratic field.  We point out that idea analogy between fundamental domain algorithms for Fuchsian groups and continued fraction algorithms goes back at leastto Eichler~\cite{Eichler}.  It seems likely therefore that many other results which follow from the theory of modular symbols should hold in the context of Shimura curves and Hilbert modular varieties as well.

\section{Cubic fields}

In this section, we tabulate some examples computed with the algorithms illustrated above.  We perform our calculations in \textsf{Magma} \cite{Magma}.  

\subsection*{First example}

We first consider the cubic field $F$ with $d_F=1101 = 3 \cdot 367$ and primitive element $w$ satisfying $w^3 - w^2 - 9w + 12=0$.  The field $F$ has Galois group $S_3$ and strict class number $1$.  The Shimura curve $X(1)=X_0^{(1)}(1)$ associated to $F$ has signature $(1; 2^2, 3^5)$---according to tables of the second author \cite{VoightShim}, this is the cubic field of strict class number $1$ with smallest discriminant such that the corresponding Shimura curve has genus $\geq 1$.  

We first compute the representative quaternion algebra $B=\quat{-1,-w^2 + w + 1}{F}$ with discriminant $\fD=(1)$ and a maximal order $\cO$, generated as a $\Z_F$-algebra by $\alpha$ and the element
\[ \frac{1}{4}\bigl( (-8w + 14) + (-2w + 4)\alpha + (-w + 2)\beta\bigr). \]
We then compute a fundamental domain for $\Gamma(1)=\Gamma_0^{(1)}(1)$; it is shown in Figure 7.1.

\begin{figure*}[h]
\begin{center}
\psset{unit=1.5in}
\begin{pspicture}(-1,-1)(1,1)
\pscircle[fillstyle=solid,fillcolor=lightgray](0,0){1}

\psclip{\pscircle(0,0){1}} \pscircle[fillstyle=solid,fillcolor=white](9.52632,0.000000){9.47368} \endpsclip
\psclip{\pscircle(0,0){1}} \pscircle[fillstyle=solid,fillcolor=white](0.00820339,1.01873){0.194613} \endpsclip
\psclip{\pscircle(0,0){1}} \pscircle[fillstyle=solid,fillcolor=white](-0.210930,0.983315){0.106768} \endpsclip
\psclip{\pscircle(0,0){1}} \pscircle[fillstyle=solid,fillcolor=white](-0.313788,0.954688){0.0994620} \endpsclip
\psclip{\pscircle(0,0){1}} \pscircle[fillstyle=solid,fillcolor=white](-0.428320,0.905832){0.0631689} \endpsclip
\psclip{\pscircle(0,0){1}} \pscircle[fillstyle=solid,fillcolor=white](-0.484547,0.877043){0.0631689} \endpsclip
\psclip{\pscircle(0,0){1}} \pscircle[fillstyle=solid,fillcolor=white](-0.511980,0.860477){0.0504584} \endpsclip
\psclip{\pscircle(0,0){1}} \pscircle[fillstyle=solid,fillcolor=white](-0.733425,0.721760){0.242590} \endpsclip
\psclip{\pscircle(0,0){1}} \pscircle[fillstyle=solid,fillcolor=white](-0.929955,0.410020){0.181472} \endpsclip
\psclip{\pscircle(0,0){1}} \pscircle[fillstyle=solid,fillcolor=white](-0.983015,0.194864){0.0654956} \endpsclip
\psclip{\pscircle(0,0){1}} \pscircle[fillstyle=solid,fillcolor=white](-0.993428,0.125088){0.0504584} \endpsclip
\psclip{\pscircle(0,0){1}} \pscircle[fillstyle=solid,fillcolor=white](-1.00387,0.0463281){0.0994620} \endpsclip
\psclip{\pscircle(0,0){1}} \pscircle[fillstyle=solid,fillcolor=white](-1.00387,-0.0463281){0.0994620} \endpsclip
\psclip{\pscircle(0,0){1}} \pscircle[fillstyle=solid,fillcolor=white](-0.993428,-0.125088){0.0504584} \endpsclip
\psclip{\pscircle(0,0){1}} \pscircle[fillstyle=solid,fillcolor=white](-0.983015,-0.194864){0.0654956} \endpsclip
\psclip{\pscircle(0,0){1}} \pscircle[fillstyle=solid,fillcolor=white](-0.929955,-0.410020){0.181472} \endpsclip
\psclip{\pscircle(0,0){1}} \pscircle[fillstyle=solid,fillcolor=white](-0.733425,-0.721760){0.242590} \endpsclip
\psclip{\pscircle(0,0){1}} \pscircle[fillstyle=solid,fillcolor=white](-0.511980,-0.860477){0.0504584} \endpsclip
\psclip{\pscircle(0,0){1}} \pscircle[fillstyle=solid,fillcolor=white](-0.484547,-0.877043){0.0631689} \endpsclip
\psclip{\pscircle(0,0){1}} \pscircle[fillstyle=solid,fillcolor=white](-0.428320,-0.905832){0.0631689} \endpsclip
\psclip{\pscircle(0,0){1}} \pscircle[fillstyle=solid,fillcolor=white](-0.313788,-0.954688){0.0994620} \endpsclip
\psclip{\pscircle(0,0){1}} \pscircle[fillstyle=solid,fillcolor=white](-0.210930,-0.983315){0.106768} \endpsclip
\psclip{\pscircle(0,0){1}} \pscircle[fillstyle=solid,fillcolor=white](0.00820339,-1.01873){0.194613} \endpsclip

\psclip{\pscircle(0,0){1}} \pscircle(9.52632,0.000000){9.47368} \endpsclip
\psclip{\pscircle(0,0){1}} \pscircle(0.00820339,1.01873){0.194613} \endpsclip
\psclip{\pscircle(0,0){1}} \pscircle(-0.210930,0.983315){0.106768} \endpsclip
\psclip{\pscircle(0,0){1}} \pscircle(-0.313788,0.954688){0.0994620} \endpsclip
\psclip{\pscircle(0,0){1}} \pscircle(-0.428320,0.905832){0.0631689} \endpsclip
\psclip{\pscircle(0,0){1}} \pscircle(-0.484547,0.877043){0.0631689} \endpsclip
\psclip{\pscircle(0,0){1}} \pscircle(-0.511980,0.860477){0.0504584} \endpsclip
\psclip{\pscircle(0,0){1}} \pscircle(-0.733425,0.721760){0.242590} \endpsclip
\psclip{\pscircle(0,0){1}} \pscircle(-0.929955,0.410020){0.181472} \endpsclip
\psclip{\pscircle(0,0){1}} \pscircle(-0.983015,0.194864){0.0654956} \endpsclip
\psclip{\pscircle(0,0){1}} \pscircle(-0.993428,0.125088){0.0504584} \endpsclip
\psclip{\pscircle(0,0){1}} \pscircle(-1.00387,0.0463281){0.0994620} \endpsclip
\psclip{\pscircle(0,0){1}} \pscircle(-1.00387,-0.0463281){0.0994620} \endpsclip
\psclip{\pscircle(0,0){1}} \pscircle(-0.993428,-0.125088){0.0504584} \endpsclip
\psclip{\pscircle(0,0){1}} \pscircle(-0.983015,-0.194864){0.0654956} \endpsclip
\psclip{\pscircle(0,0){1}} \pscircle(-0.929955,-0.410020){0.181472} \endpsclip
\psclip{\pscircle(0,0){1}} \pscircle(-0.733425,-0.721760){0.242590} \endpsclip
\psclip{\pscircle(0,0){1}} \pscircle(-0.511980,-0.860477){0.0504584} \endpsclip
\psclip{\pscircle(0,0){1}} \pscircle(-0.484547,-0.877043){0.0631689} \endpsclip
\psclip{\pscircle(0,0){1}} \pscircle(-0.428320,-0.905832){0.0631689} \endpsclip
\psclip{\pscircle(0,0){1}} \pscircle(-0.313788,-0.954688){0.0994620} \endpsclip
\psclip{\pscircle(0,0){1}} \pscircle(-0.210930,-0.983315){0.106768} \endpsclip
\psclip{\pscircle(0,0){1}} \pscircle(0.00820339,-1.01873){0.194613} \endpsclip
\pscircle(0,0){1}
\end{pspicture} \\
\textbf{Figure 7.1}: A fundamental domain for $X(1)$ for $d_F=1101$
\end{center}
\end{figure*}

We may take $\mu=\beta$ as an element to represent complex conjugation; indeed, $\beta^2=-w^2 + w + 1 \in \Z_F^*$ is a unit.  Since the spaces $H^1(\Gamma,\Z)^{\pm}$ are one-dimensional, the Hecke operators $T_\fp$ act by scalar multiplication, and the eigenvalues are listed in Table 7.2: for each prime $\fp$ of $\Z_F$ with $N(\fp) < 50$, we list a generator $\pi$ of $\fp$ and the eigenvalue $a(\fp)$ of $T_\fp$.  

\begin{table}[h]
\[ \begin{array}{cc|cc}
N \fp & \pi & a(\fp) & \#E(\FF_{\fp}) \\
\hline
2\rule{0pt}{2.5ex}   & w - 2 & 0 & 3 \\
3 & w - 3 & -3 & 7 \\
3 & w - 1 & -1 & 5 \\
4 & w^2 + w - 7 & -3 & 8 \\
19 & w + 1 & -6 & 26 \\
23 & w^2 - 2w - 1 & 6 & 18 \\
31 & 2w^2 - 19 & 3 & 29 \\
31 & w^2 - 5 & 0 & 32 \\
31 & 3w - 5 & 4 & 28 \\
41 & w^2 + 2w - 7 & 0 & 42 \\
43 & w^2 - 11 & 9 & 35 \\
47 & 3w - 7 & -9 & 57
\end{array} \] 
\textbf{Table 7.2}: Hecke eigenvalues for the group $\Gamma(1)$ for $d_F=1101$
\end{table}

The curve $X(1)$ has modular Jacobian $E$ and we have $\# E(\FF_\fp)=N \fp +1-a(\fp)$, so we list these values in Table 7.2 as well.  By the reduction theory of Carayol \cite{Carayol}, we know that $E$ is an elliptic curve over $F$ with everywhere good reduction.  (Compare this result with calculations of Demb\'el\'e and Donnelly \cite{DD}.)  We note that since the $a(\fp)$ for the primes $\fp$ with $N\fp=31$ are not equal, the curve $E$ does not arise as the base change of a curve defined over $\Q$. 

To find a candidate curve $E$, we begin by searching for a curve over $F$ with everywhere good reduction.  We follow the methods of Cremona and Lingham \cite{CremonaLingham}.  

\begin{remark}
One possible alternative approach to find an equation for the curve $E$ would be to use the data computed above to give congruence conditions on a minimal Weierstrass model for $E$
\[ E:y^2 + a_1 xy + a_3 y = x^3 + a_2 x^2 + a_4 x^4 + a_6 \] 
with $a_i \in \Z_F$.  For example, by coordinate change for $y$, without loss of generality we may assume that $a_1,a_3$ are chosen from representatives of the set $\Z_F/2\Z_F$, and since the reduction of $E$ modulo the prime $\fp$ with $N\fp=2$ is supersingular, and a Weierstrass model for such curve is of the form $y^2+y=f(x)$, we may assume that $a_1 \equiv 0 \pmod{w-2}$ and $a_3 \equiv 1 \pmod{w-2}$, leaving $4$ possibilities each for $a_1,a_3$.  In a similar way, we obtain further congruences.  In our case, this approach fails to find a model, but we expect it will be useful in many cases.
\end{remark}

\begin{remark}
When $F$ is a real quadratic field, there are methods of Demb\'el\'e which apply \cite{DembeleRealQuadratic} to find an equation for the curve $E$ by computing the real periods of $E$.
\end{remark}

In our situation, where $F$ has (strict) class number $1$, we conclude (see Cremona and Lingman \cite[Propositions 3.2--3.3]{CremonaLingham}) that there exists $w \in \Z_F^*/\Z_F^{*6}$ and an integral point on the elliptic curve $E_w:y^2=x^3-1728w$ such that $E$ has $j$-invariant $j=x^3/w=1728+y^2/w$.  In our situation, we need not provably enumerate all such curves and content ourselves to find as many such integral points as we can.  Indeed, we find for the unit $w=-1506w^2 + 6150w - 5651 \in \Z_F^*$ that the curve $E_w$ has rank $3$ and we find an integral point
\[ (-11w^2 - 24w + 144, -445w^2 + 1245w - 132) \in E_w(F); \]
this point corresponds to a curve with $j$-invariant 
\[ (-2w^2-4w+7)^9(2w^2+w-17)^3(w-2)^{18}(w-3)^6 = -1805w^2 - 867w + 14820 \]
where we note that the first two terms are units and recall that $N(w-2)=2$, $N(w-3)=-3$.  We then find an appropriate quadratic twist of this curve which has conductor $(1)=\Z_F$ as follows:
\[ A: y^2 + w(w + 1)xy + (w + 1)y = x^3 + w^2x^2 + a_4x + a_6 \]
where
\begin{center}
\footnotesize{$a_4 = -139671409350296864w^2 - 235681481839938468w + 623672370161912822$}
\end{center}
and
\begin{center}
\footnotesize{$a_6 = 110726054056401930182106463w^2 + 186839095087977344668356726w - 494423184252818697135532743$}.
\end{center}
We verify that $\#A(\F_\fp)=N(\fp)+1-a(\fp)$ in agreement with the above table, so this strongly suggests that the Jacobian $E=J(1)$ of $X(1)$ is isogenous to $A$ (and probably isomorphic to $A$).

To prove that in fact $E$ is isogeneous to $A$, we use the method of Faltings and Serre.  (For an exposition of this method, see Sch\"utt \cite[\S 5]{Schutt} and Dieulefait, Guerberoff, and Pacetti \cite[\S 4]{Dieule}, and the references contained therein.)  

Let $\rho_E,\rho_A: \Gal(\overline{F}/F) \to \GL_2(\Z_2)$ be the $2$-adic Galois representations associated to the $2$-adic Tate modules of $E$ and $A$, respectively, and let $\overline{\rho}_E,\overline{\rho}_A : \Gal(\overline{F}/F) \to \GL_2(\F_2)$ be their reductions modulo $2$.  We will show that we have an isomorphism $\overline{\rho}_E \cong \overline{\rho}_A$ of absolutely irreducible representations and then lift this to an isomorphism $\rho_E \cong \rho_A$ by comparing traces of $\rho_E$ and $\rho_A$.  It then follows from work of Faltings that $E$ is isogeneous to $A$.

First we show that the representations $\overline{\rho}_E$ and $\overline{\rho}_A$ are absolutely irreducible and isomorphic.  (It is automatic that they have the same determinant, the cyclotomic character.)  For $E$, it is clear from Table 7.2 that the image of $\overline{\rho}_E$ contains elements both of even and odd order and so must be all of $\GL_2(\F_2) \cong S_3$.  For $A$, we verify that the $2$-division polynomial is irreducible, and adjoining a root of this polynomial to $F$ we obtain a field which is more simply given by the polynomial
\[ x^3 + (w + 1)x^2 + (w^2 + w - 5)x + (w^2 + w - 7) \]
and with relative discriminant $4\Z_F$; the splitting field $L$ of this polynomial indeed has Galois group isomorphic to $S_3$.

At the same time, the field cut out by $\overline{\rho}_E$ is an $S_3$-extension of $F$ which is unramified away from $2$.  We may enumerate all such fields using class field theory, as follows.  Any such extension has a unique quadratic subextension which is also unramified away from $2$ and by Kummer theory is given by adjoining the square root of a $2$-unit $\alpha \in \Z_F$.  Since $2\Z_F=\fp_2\fp_2'$ with $\fp_2=(w-2)\Z_F$, $\fp_2'=(w^2+w-7)\Z_F$ (of inertial degrees $1,2$, respectively), we see that there are $31=2^{3+2}-1$ possibilities for $\alpha$.  Now for each quadratic extension $K$, we look for a cyclic cubic extension which is unramified away from $2$ and which generates a Galois field over $F$.  By class field theory, we computing the complete list, and in fact we find a unique such field and thereby recover the field $L$.  The field $L$ arises from the quadratic subfield $K=F(\sqrt{u})$, where $u=-19w^2-32w+85 \in \Z_F^*$; the class group of $K$ is trivial, but the ray class group modulo $\fp_2=(w-2)\Z_F$ is isomorphic to $\Z/3\Z$.  We have therefore indeed shown that $\overline{\rho}_E$ and $\overline{\rho}_A$ are isomorphic.

Next, we lift this isomorphism to one between $\rho_E$ and $\rho_A$.  Following Faltings and Serre, we must first compute all quadratic extensions $M$ of $L$ which are unramified away from $2$.  We find that the $2$-unit group of $\Z_L$ modulo squares has rank $16$ over $\Z/2\Z$; in other words, we have $16$ independent quadratic characters of $L$.  To each odd ideal $\fA$ of $\Z_L$, we associate the corresponding vector $v(\fA) = (\chi(\fA))_{\chi} \in \{\pm 1\}^{16}$ of values of these $16$ characters.  Then for each $v \in \{\pm 1\}^{16}$, we need to find an odd ideal $\fA$ of $\Z_L$ such that $v(\fA)=v$ and then verify that $\rho_A(\Frob_\fa)=\rho_E(\Frob_\fa)$, where $\fa = \Z_F \cap \fA$.  For this, it is enough to verify that $\tr \rho_A(\Frob_\fa)=\tr \rho_E(\Frob_\fa)$.  Moreover, by multiplicativity, thinking of $\{\pm 1\}^{16}$ as an $\F_2$-vector space in a natural way, it is enough to verify this for a set of odd primes $\fP$ of $\Z_L$ such that the values $v(\fP)$ span $\{\pm 1\}^{16}$.  We find that the set of primes $\fP$ of $\Z_L$ which lie above primes $\fp$ of $\Z_F$ of norm at most $239$ will suffice for this purpose, and for this set we indeed verify that the Hecke eigenvalues agree.  (In other words, we verify that $\#E(\F_{\fp})=\#A(\F_{\fp})$ for enough primes $\fp$ of $\Z_F$ to ensure that in fact equality holds for all odd primes $\fp$ of $\Z_F$, whence $\rho_J \cong \rho_A$.)  Therefore $E$ is isogeneous to $A$, and we have explicitly identified the isogeny class of the Jacobian of the Shimura curve $X(1)$ over $F$.

\subsection*{Second example}

As a second example, we consider the Galois cubic field with $d_F=1369=37^2$, so that $F$ is the (unique totally real) cubic field in $\QQ(\zeta_{37})$.  The field $F$ is generated by an element $w$ satisfying $w^3-w^2-12w-11$.  Here, the associated Shimura curve of level and discriminant $1$ has signature $(1;2^3,3^3)$.  Computing as above, we obtain the Hecke eigenvalues listed in Table 7.3.

\begin{table}[h]
\[ \begin{array}{c|cc}
N \fp & a(\fp) & \#E(\FF_p) \\
\hline
8\rule{0pt}{2.5ex} & -5 & 14 \\
11 & -2 & 14 \\
23 & -4 & 28 \\
27 & 0 & 28 \\
29 & 9 & 21 \\
31 & -10 & 42 \\
37 & -11 & 49 \\
43 & 2 & 42 \\
47 & 6 & 42
\end{array} \] 
\textbf{Table 7.3}: Hecke eigenvalues for the group $\Gamma(1)$ for $d_F=1369$
\end{table}

Let $E$ denote the Jacobian of the Shimura curve $X(1)$ under consideration.  Since $\#E(\FF_p)$ is always divisible by $7$, it is reasonable to believe that $E$ has a nontrivial $7$-torsion point.  A quick search for $F$-rational points on $X_1(7)$ (using the Tate model) yields a curve $A$ with minimal model
\[ A: y^2 + (w^2 + w + 1)xy = x^3 + (-w - 1)x^2 + (256w^2 + 850w + 641)x + (5048w^2 + 16881w + 12777) \]
with $\#A(\FF_p)$ as above.  Since $A[7](F)$ is nontrivial, an argument of Skinner-Wiles \cite{SkinnerWiles} shows that in fact $A$ is isogenous to $E$.  

We note that the eigenvalue $a(\fp)$ does not depend on the choice of prime $\fp$ above $p$, which suggests that $E$ should come as a base change from a curve defined over $\Q$.  Indeed, it can be shown by Galois descent (arising from the functoriality of the Shimura-Deligne canonical model) that the curve $X(1)$ has field of moduli equal to $\Q$.  Looking at the curves of conductor $1369$ in Cremona's tables \cite{Cremona}, we find that $A$ is in fact the base change to $F$ of the curve
\[ \text{\textsf{1369b1}}: y^2 + xy  = x^3 - x^2 + 3166x - 59359. \]

\subsection*{Third example}

We conclude with a third example for which the dimension of the space of cuspforms is greater than $1$.  The first cubic field $F$ for which this occurs has $d_F=961=31^2$, and the signature of the corresponding Shimura curve $X(1)$ is $(2;2^4,3)$.  Since this field $F$ is Galois, by Galois descent the corresponding $L$-function is a base change from $\QQ$.  We list the characteristic polynomial $\chi(T_\fp)$ of Hecke operators $T_\fp$ in Table 7.4.

\begin{table}[h]
\[ \begin{array}{c|cc}
N \fp & \chi(T_\fp) \\
\hline
2\rule{0pt}{2.5ex} & x^2+2x-1 \\
23 & x^2+8x+16 \\
27 & x^2-4x-28 \\
29 & x^2+8x+8 \\
31 & x^2 + 20x + 100 \\
47 & x^2 - 8x - 16
\end{array} \] 
\textbf{Table 7.4}: Hecke eigenvalues for the group $\Gamma(1)$ for $d_F=961$
\end{table}

We see that $\Q(T_\fp)=\Q(\sqrt{5})$, so the Jacobian $J$ of $X(1)$ is an abelian surface with real multiplication by $\Q(\sqrt{5})$.  Computing via modular symbols the space of newforms for the classical modular group $\Gamma_0(961)$, we find that these characteristic polynomials indeed match a newform with $q$-expansion
\[ q + wq^2 + wq^3 + (-2w - 1)q^4 + \dots -4q^{23} + \dots - (2w + 6)q^{29} + \dots \]
where $w^2+2w-1=0$.
It is interesting to note that there are six $2$-dimensional (irreducible) new Hecke eigenspaces for $\Gamma_0(961)$ and only one of the corresponding forms corresponds to an abelian surface which acquires everywhere good reduction over $F$.

\subsection*{Comments}

Combining the above methods with those of Demb\'el\'e and Donnelly, one can systematically enumerate Hilbert modular forms over a wide variety of totally real fields and conductors.  This project has been initiated by Donnelly and the second author \cite{VoightDonnelly}, and we refer to this upcoming work for many, many more examples for fields up to degree $6$.

We conclude with a few comments on the efficiency of our algorithms.  Unfortunately, we have no provable assessment of the running time of our method.  In practice, the computation of a fundamental domain seems to be the most (unpredictably) time-consuming step, taking on the order of minutes on a standard PC for cubic and quintic fields of small discriminant; but this step should be considered a precomputation step as it need be done only once for each field $F$.  The computation of a single Hecke operator takes on the order of seconds (for the examples above) to a few minutes (for quintic fields), the most onerous step being the principalization step (Proposition \ref{KVprinc}) and the bookkeeping involved in working with the induced module; with further careful optimization, we believe that this running time can be substantially lowered.

\end{document}